\documentclass[a4paper,12pt,twoside]{article}
\usepackage{amsmath,amssymb,ifthen}
\usepackage[amsmath,thmmarks]{ntheorem}
\usepackage{paper}
\usepackage[all]{xy}
\newcommand{\p}{\mathfrak{p}}
\newcommand{\m}{\mathfrak{m}}
\DeclareMathOperator{\Ind}{Ind}

\DeclareMathOperator{\Def}{Def}
\DeclareMathOperator{\Tor}{Tor}

\DeclareMathOperator{\rank}{rank}
\DeclareMathOperator{\cha}{char}
\newcommand{\yleftarrow}[2][]{\settowidth{\dimen0}{\ensuremath{#1}}\settowidth{\dimen1}{\ensuremath{#2}}%
\xleftarrow[\ifthenelse{\lengthtest{\dimen0>12.5pt}}{#1}{\makebox[12.5pt]{\ensuremath{\scriptstyle #1}}}]%
{\ifthenelse{\lengthtest{\dimen1>12.5pt}}{#2}{\makebox[12.5pt]{\ensuremath{\scriptstyle #2}}}}}
\makeatletter
\newtheoremstyle{nonumberplainwithoutbrackets}%
  {\item[\theorem@headerfont\hskip\labelsep ##1\theorem@separator]}%
  {\item[\theorem@headerfont\hskip \labelsep ##1\ ##3\theorem@separator]}
\makeatother
\theoremstyle{nonumberplainwithoutbrackets}
\newtheorem{prfwithname}{Proof}

\title{Non-cuspidality outside the middle degree of $\ell$-adic cohomology of the Lubin-Tate tower}
\author{Yoichi Mieda}

\begin{document}

\maketitle

\begin{firstfootnote}
 Graduate School of Mathematics, Kyushu University, 6--10--1 Hakozaki, Higashi-ku, Fukuoka, 812--8581 Japan

 E-mail address: \texttt{mieda@math.kyushu-u.ac.jp}

\end{firstfootnote}

\begin{abstract}
 In this article, we consider the representations of the general linear group over a non-archimedean local field
 obtained from the vanishing cycle cohomology of the Lubin-Tate tower. We give an easy and direct proof
 of the fact that no supercuspidal representation appears as a subquotient of such representations
 unless they are obtained from the cohomology of the middle degree.
 Our proof is purely local and does not require Shimura varieties.
\end{abstract}

\section{Introduction}
Let $F$ be a non-archimedean local field, i.e., a complete discrete valuation field with finite residue field,
and $\mathcal{O}$ the ring of integers of $F$. For an integer $d\ge 1$,
we consider the universal deformation space over $(\mathcal{O}^{\mathrm{ur}})^{\wedge}$, 
the completion of the maximal unramified extension of $\mathcal{O}$, of formal $\mathcal{O}$-modules of height $d$.
It is called the Lubin-Tate space.
By adding Drinfeld level structures, we get a tower over it, which is called the Lubin-Tate tower.

Let $k$ be a non-negative integer. Using the $k$th vanishing cycle cohomology of the Lubin-Tate tower,
we obtain a representation of $\GL_d(F)\times D^\times\times W_F$, which we denote by $H^k_{\mathrm{LT}}$
(\cf \cite{Carayol}; see also Definition \ref{defn:H_LT}, where we give the construction of
 $H^k_{\mathrm{LT}}$ as a $\GL_d(F)$-representation).
Here $D$ is the central division algebra over $F$ with invariant $1/d$, and $W_F$ is the Weil group of $F$.
The representation $H^k_{\mathrm{LT}}$ is very important for the study of the local Langlands correspondence;
in fact, it is known that for a supercuspidal representation $\pi$ of $\GL_d(F)$ the $\pi$-isotypic component
of $H^{d-1}_{\mathrm{LT}}$ can be described by the local Langlands correspondence and 
the local Jacquet-Langlands correspondence. This was known as non-abelian Lubin-Tate theory
or the conjecture of Deligne-Carayol.

The aim of this paper is to give an easy and direct proof of the following theorem:

\begin{thm}[Theorem \ref{thm:main-thm}]\label{thm:intro}
 For $k\neq d-1$, the $\GL_d(F)$-representation $H^k_{\mathrm{LT}}$ 
 has no supercuspidal representation as a subquotient.
\end{thm}

Actually, the theorem has already been proved;
for the case $\cha F\neq 0$ it was due to P.~Boyer \cite[Th\'eor\`eme 3.2.4]{Boyer},
and for the case $\cha F=0$ it was remarked by G.~Faltings \cite{Faltings}. See also \cite{Boyer2, Boyer3}, where
the $\GL_d(F)\times D^\times\times W_F$-representation $H^k_{\mathrm{LT}}$ is completely described for every $k$.
These proofs are accomplished by relating the Lubin-Tate tower to an integral model of a certain Shimura variety
or the moduli space of $\mathcal{D}$-elliptic sheaves with level structures; hence they are somewhat indirect.
It seems natural to look for a more direct proof. In this direction, there are some previous works;
in \cite{Yoshida}, T.~Yoshida considered the ``depth 0'' case (\cite[Theorem 6.16 (ii)]{Yoshida});
in \cite{Strauch}, under some finiteness assumptions (\cite[4.2.7 (H)]{Strauch}),
M.~Strauch gave a proof of Theorem \ref{thm:intro} (\cite[Section 4.3]{Strauch}).
Our method is inspired by the paper \cite{Strauch}, but works unconditionally.
The point is that we remain in the category of schemes, that is, we use neither formal schemes nor rigid spaces.
We use Gabber's new result (\cite{Gabber}) in order to obtain the finiteness result needed for our proof.

We sketch the outline of this paper. In Section 2, we gather some general results on the nearby cycle
cohomology of a complete local ring over a strict local discrete valuation ring. In Section 3, we recall
definitions and notation concerning the Lubin-Tate tower, and state our main result. In Section 4,
we give a proof of our main theorem.

\medbreak
\noindent\textbf{Acknowledgment}\quad
The author would like to thank Tetsushi Ito and Matthias Strauch 
for reading a manuscript and giving helpful comments.
He was supported by the Japan Society for the Promotion of Science Research Fellowships for
Young Scientists.

\medbreak
\noindent\textbf{Convention}\quad Every sheaf and cohomology are considered on the \'etale site of a scheme.

\section{Complements on \'etale cohomology}
Let $V$ be a strict local excellent discrete valuation ring and put $S=\Spec V$.
 We denote the closed (resp.\ generic) point by $s$ (resp.\ $\eta$). 
For a scheme $X$ over $S$, we write $X_s$ for its special fiber and $X_\eta$ for its generic fiber.
Fix a prime number $\ell$ which is invertible in $V$ and a positive integer $n$.
Put $\Lambda=\Z/\ell^n\Z$.

Let $\mathcal{X}$ be a scheme separated of finite type over $S$. For a closed point $x$ of $\mathcal{X}$
lying over $s$, denote by $A=\mathcal{O}_{\mathcal{X},x}^\wedge$
the completion of the local ring $\mathcal{O}_{\mathcal{X},x}$ and
put $X=\Spec A=\Spec \mathcal{O}_{\mathcal{X},x}^\wedge$.

\begin{prop}\label{prop:regular-base-change}
 For every integer $k\ge 0$, the nearby cycle sheaf $R^k\psi_X\Lambda$ is a constructible sheaf on $X_s$.
\end{prop}

\begin{prf}
 Consider the following commutative diagram whose rectangles are cartesian:
 \[
  \xymatrix{%
 X_s\ar[r]^-{\hat{\imath}}\ar[d]^-{\varphi_s}& X\ar[d]^-{\varphi}& X_{\overline{\eta}}\ar[l]_-{\hat{\jmath}}\ar[d]^-{\varphi_{\overline{\eta}}}\\
 \mathcal{X}_s\ar[r]^-i& \mathcal{X}& \mathcal{X}_{\overline{\eta}}\ar[l]_-{j}\lefteqn{.}
 }
 \]
 Here $\overline{\eta}$ is a geometric point lying over $\eta$.
 Then by definition we have $R\psi_{\mathcal{X}}\Lambda=i^*Rj_*\Lambda$ and
 $R\psi_X\Lambda=\hat{\imath}^*R\hat{\jmath}_*\Lambda$. Since $\mathcal{X}$ is an excellent scheme,
 $\varphi$ is a regular morphism. Therefore, by the regular base change theorem 
 (\cite[Corollary 7.1.6]{Fujiwara}, see also \cite[Proposition 4.5]{Riou}), we have
 \[
 \varphi_s^*R\psi_{\mathcal{X}}\Lambda=\varphi_s^*i^*Rj_*\Lambda=\hat{\imath}^*\varphi^*Rj_*\Lambda
 \cong \hat{\imath}^*R\hat{\jmath}_*\Lambda=R\psi_X\Lambda.
 \]
 Hence we have $R^k\psi_{\mathcal{X}}\Lambda\cong \varphi^*R^k\psi_X\Lambda$.
 Since $R^k\psi_{\mathcal{X}}\Lambda$ is a constructible sheaf on $\mathcal{X}_s$
 (\cite[Finitude, Th\'eor\`eme 3.2]{SGA4+1/2}), so is $R^k\psi_X\Lambda$
 (\cite[Expos\'e IX, Proposision 2.4 (iii)]{SGA4}).
\end{prf}

\begin{prop}\label{prop:finiteness-local-coh}
 Let $Z$ be a locally closed subscheme of $X_s$. 
 Then for every integer $k$, the local cohomology $H^k_Z(X_s,R\psi_X\Lambda)$
 is a finitely generated $\Lambda$-module.
\end{prop}

\begin{prf}
 By Proposition \ref{prop:regular-base-change}, it suffices to show that for every constructible sheaf
 $\mathcal{F}$ on $X_s$, $H^k_Z(X_s,\mathcal{F})$ is a finitely generated $\Lambda$-module.
 First we consider the case where $Z$ is open in $X_s$.
 Denote by $j\colon Z\hooklongrightarrow X_s$ the natural open immersion.
 Since $X_s$ is noetherian, $j$ is quasi-compact, and thus is of finite type. Therefore by Gabber's finiteness
 theorem (\cite{Gabber}), $R^kj_*j^*\mathcal{F}$ is a constructible sheaf on $X_s$.
 Thus $H^k_Z(X_s,\mathcal{F})=H^k(Z,j^*\mathcal{F})=H^k(X_s,Rj_*j^*\mathcal{F})=(R^kj_*j^*\mathcal{F})_x$
 is a finitely generated $\Lambda$-module, as desired. Next we consider the case where $Z$ is closed in $X_s$.
 By the exact sequence $H^{k-1}(X_s\setminus Z,\mathcal{F})\longrightarrow H^k_Z(X_s,\mathcal{F})\longrightarrow H^k(X_s,\mathcal{F})$, we may reduce to the open case. Finally we consider the general case. Let $\overline{Z}$ be
 the closure of $Z$ in $X_s$. Then $Z$ is open in $\overline{Z}$.
 By the exact sequence $H^k_{\overline{Z}}(X_s,\mathcal{F})\longrightarrow H^k_Z(X_s,\mathcal{F})\longrightarrow H^{k+1}_{\overline{Z}\setminus Z}(X_s,\mathcal{F})$, we may reduce to the closed case.
 Now the proof is complete.
\end{prf}

\begin{rem}\label{rem:without-Gabber}
 If $Z$ comes from a locally closed subscheme of $\mathcal{X}_s$, then by the similar method as in the proof
 of Proposition \ref{prop:regular-base-change}, we may avoid Gabber's deep theorem.
\end{rem}

\begin{prop}\label{prop:coh-dim}
 Assume that $\mathcal{X}_\eta$ is smooth of pure dimension $d-1$. Then we have $(R^k\psi_X\Lambda)_x=0$ for
 $k\ge d$ and $H^k_x(X_s,R\psi_X\Lambda)=0$ for $k\le d-2$.
\end{prop}

\begin{prf}
 We use the same notation as in the proof of Proposition \ref{prop:regular-base-change}.
 First, \cite[Expos\'e XIII, 2.1.13]{SGA7} ensures the vanishing of $(R^k\psi_X\Lambda)_x\cong (R^k\psi_{\mathcal{X}}\Lambda)_x$ for $k\ge d$.

 Denote the closed immersion $x\hooklongrightarrow \mathcal{X}_s$ (resp.\ $x\hooklongrightarrow X_s$)
 by $i_x$ (resp. $\hat{\imath}_x$).
 Since $\varphi_s$ is a regular morphism, we have 
 $\hat{\imath}_x^!R\psi_X\Lambda\cong \hat{\imath}_x^!\varphi_s^*R\psi_{\mathcal{X}}\Lambda\cong i_x^!R\psi_{\mathcal{X}}\Lambda$ (\cite[Corollaire 4.7]{Riou}). In other words, we have
 $H^k_x(X_s,R\psi_X\Lambda)\cong H^k_x(\mathcal{X}_s,R\psi_{\mathcal{X}}\Lambda)$.
 On the other hand, since we have
 \[
  D_x(i_x^!R\psi_{\mathcal{X}}\Lambda)\cong i_x^*D_{\mathcal{X}_s}(R\psi_{\mathcal{X}}\Lambda)
 \cong i_x^*R\psi_{\mathcal{X}}D_{\mathcal{X}_\eta}(\Lambda)\cong 
 i_x^*R\psi_{\mathcal{X}}\Lambda(d-1)[2d-2]
 \]
 (for the second isomorphy, see \cite[Th\'eor\`eme 4.2]{Illusie}), the two $\Lambda$-modules
 $H^k_x(\mathcal{X}_s,R\psi_{\mathcal{X}}\Lambda)$ and 
 $(R^{2d-2-k}\psi_{\mathcal{X}}\Lambda)_x(d-1)$ are dual to each other.
 Therefore $H^k_x(\mathcal{X}_s,R\psi_{\mathcal{X}}\Lambda)=0$ for $2d-2-k\ge d$, i.e., $k\le d-2$.
 Now the proof is complete.
\end{prf}

Let $B$ be a $V$-algebra and $G$ a finite group of $V$-algebra isomorphisms of $B$.
Put $Y=\Spec B$ and $Y'=\Spec B^G$. Denote the natural morphism $Y\longrightarrow Y'$ by $\pi$.
Let $W'$ be a locally closed subscheme of $Y'_s$ and put $W=\pi^{-1}(W')$.

\begin{prop}\label{prop:invariant}
 Assume that the cardinality of $G$ is prime to $\ell$, $\pi$ is finite, and 
 $\pi_\eta\colon Y_\eta\longrightarrow Y_\eta'$ is \'etale.
 Then for every integer $k$ we have an isomorphism 
 \[
  H^k_{W'}(Y_s',R\psi_{Y'}\Lambda)\cong H^k_W(Y_s,R\psi_{Y}\Lambda)^G.
 \]
\end{prop}

\begin{prf}
 For a $\Lambda[G]$-sheaf $\mathcal{F}$ over a scheme, we denote by $\mathcal{F}^G$
 the $\Lambda$-sheaf defined by $U\longmapsto \mathcal{F}(U)^G$. 
 It is easy to see that this functor is exact (by the assumption on the cardinality of $G$)
 and commutes with pull-back and derived push-forward. Therefore it also commutes with
 local cohomology and nearby cycle functor.

 Since $\pi_{\eta}\colon Y_{\eta}\longrightarrow Y'_{\eta}$ is
 an \'etale Galois covering with the Galois group $G$, we have
 $(\pi_{\eta*}\Lambda)^G\cong \Lambda$.
 Therefore we have
 \begin{align*}
 H^k_W(Y_s,R\psi_Y\Lambda)^G&\cong H^k_{W'}\bigl(Y'_s,\pi_{s*}R\psi_Y\Lambda\bigr)^G\cong 
 H^k_{W'}\bigl(Y'_s,R\psi_Y\pi_{\eta*}\Lambda\bigr)^G\\
 &\cong H^k_{W'}\bigl(Y'_s,R\psi_{Y'}(\pi_{\eta*}\Lambda)^G\bigr)\cong H^k_{W'}\bigl(Y'_s,R\psi_{Y'}\Lambda\bigr),
 \end{align*}
 as desired.
\end{prf}

\begin{rem}
 By the standard method, we may prove the analogous results as Proposition \ref{prop:regular-base-change}, 
 \ref{prop:finiteness-local-coh}, \ref{prop:coh-dim}, \ref{prop:invariant} for $\Lambda=\Q_\ell$ or 
 $\Lambda=\overline{\Q}_\ell$, an algebraic closure of $\Q_\ell$.
\end{rem}

\section{Review of the Lubin-Tate tower}
In this section, we briefly recall basic definitions about the Lubin-Tate tower.
First we will introduce notation used in the sequel. Let $F$ be a complete discrete valuation field
with finite residue field $\F_q$. We denote the ring of integers of $F$ by $\mathcal{O}$ and fix a uniformizer
$\varpi$ of $\mathcal{O}$. 
Let $V=(\mathcal{O}^{\mathrm{ur}})^\wedge$ be the completion of the maximal unramified extension of $\mathcal{O}$.
Denote the residue field of $V$ by $\F$. As in the previous section, we write $s$ (resp.\ $\eta$) for
the closed (resp.\ generic) point of $S=\Spec V$.

Fix an integer $d\ge 1$. Let $\GL_d(F)^0$ (resp.\ $\GL_d(F)^{00}$) be the subgroup of $\GL_d(F)$
consisting of elements $g\in \GL_d(F)$ such that the normalized valuation of $\det g$ is divisible by $d$
(resp.\ equal to $0$).
Put $K_0=\GL_d(\mathcal{O})$.
For an integer $m\ge 1$, we put $K_m=\Ker (\GL_d(\mathcal{O})\longrightarrow \GL_d(\mathcal{O}/\varpi^m\mathcal{O}))$.
For $m\ge 0$, $K_m$ is a compact open subgroup of $\GL_d(F)$ which is contained in $\GL_d(F)^{00}$.

Let $\mathbb{G}$ be a one-dimensional formal $\mathcal{O}$-module over $\F$ with $\mathcal{O}$-height $d$.
It is known to be unique up to isomorphism. 
Let $\mathcal{C}$ be the category whose objects are complete noetherian local $V$-algebras with
residue fields isomorphic to $\F$ by the structure morphisms, and whose morphisms are local $V$-algebra
homomorphisms. For an integer $m\ge 0$, consider the covariant functor 
$\Def_m\colon \mathcal{C}\longrightarrow \mathbf{Set}$ that associates $A$ with the set of
isomorphism classes of triples $(G,\iota,\phi)$,
where $G$ is a one-dimensional formal $\mathcal{O}$-module over $A$, 
$\iota\colon \mathbb{G}\yrightarrow{\cong}G\otimes_A\F$ is an isomorphism of formal $\mathcal{O}$-modules over $\F$,
and $\phi\colon (\varpi^{-m}\mathcal{O}/\mathcal{O})^d\longrightarrow \m_A$ (here we fix a coordinate of $G$ and identify
$G(A)$ with $\m_A$, the maximal ideal of $A$)
is a Drinfeld level $m$ structure of $G$.
Recall that a Drinfeld level $m$ structure of $G$ is an $\mathcal{O}$-homomorphism from $(\varpi^{-m}\mathcal{O}/\mathcal{O})^d$ to
$\m_A$ endowed with an $\mathcal{O}$-module structure by $G$ such that
\[
 \prod_{a\in (\varpi^{-m}\mathcal{O}/\mathcal{O})^d}\bigl(T-\phi(a)\bigr)=[\varpi^m](T)g(T)
\]
for some $g(T)\in A[[T]]^\times$. This is equivalent to saying that the element $[\varpi](T)$ is  
divisible by $\prod_{a\in (\varpi^{-1}\mathcal{O}/\mathcal{O})^d}(T-\phi(a))$ (\cite[Corollary II.2.3]{Harris-Taylor}).
We define $(G,\iota,\phi)\cong (G',\iota',\phi')$ in the obvious manner.
For $m\le m'$, we have a natural morphism of functors $\Def_{m'}\longrightarrow \Def_m$
by restricting level $m'$ structures to $(\varpi^{-m}\mathcal{O}/\mathcal{O})^d\subset (\varpi^{-m'}\mathcal{O}/\mathcal{O})^d$.

The following fundamental theorem is due to Lubin, Tate and Drinfeld (\cite[Proposition 4.3, 1), 2)]{Drinfeld}):

\begin{thm}
 \begin{enumerate}
  \item For every $m\ge 0$, the functor $\Def_m$ is represented by a $d$-dimensional
	regular local $V$-algebra $A_m$.
	Moreover, as a local $V$-algebra, $A_0$ is isomorphic to $V[[T_1,\ldots,T_{d-1}]]$.
  \item Let $G^{\mathrm{univ}}$ be the universal formal $\mathcal{O}$-module over $A_0$ and
	$\phi_m^{\mathrm{univ}}\colon (\varpi^{-m}\mathcal{O}/\mathcal{O})^d\longrightarrow \m_{A_m}$ the universal
	Drinfeld level $m$ structure of $G^{\mathrm{univ}}\otimes_{A_0}A_m$.
	Then, for $m\ge 1$ and an $\mathcal{O}/\varpi^m\mathcal{O}$-basis $e_1,\ldots,e_d$ of $(\varpi^{-m}\mathcal{O}/\mathcal{O})^d$, 
	$\phi^{\mathrm{univ}}_m(e_1),\ldots,\phi^{\mathrm{univ}}_m(e_d)\in \m_{A_m}$ form a system of 
	regular	parameters of $A_m$.
 \end{enumerate}
\end{thm}

Put $X_m=\Spec A_m$ and denote the unique closed point of $X_m$ by $x_m$.
The morphism of functors $\Def_{m'}\longrightarrow \Def_m$ induces a morphism
$p_{mm'}\colon X_{m'}\longrightarrow X_m$ over $S=\Spec V$. It is finite and flat 
(\cite[Proposition 4.3, 3)]{Drinfeld}).
The $S$-scheme $X_0$ is sometimes called the Lubin-Tate space and
the projective system of schemes $(X_m)_{m\ge 0}$ is called the Lubin-Tate tower.

In fact, $A_m$ is algebraizable, that is, it can be obtained from a scheme of finite type over $V$
by taking completion at a closed point on the special fiber.

\begin{prop}\label{prop:algebraization-A_m}
 For every integer $m\ge 0$,
 there exist a scheme $\mathfrak{M}_m$ of finite type over $V$
 with purely $(d-1)$-dimensional smooth generic fiber, a closed point $y$ of $\mathfrak{M}_{m,s}$, 
 and a local $V$-algebra isomorphism $A_m\cong \mathcal{O}_{\mathfrak{M}_m,y}^{\wedge}$, where ${}^\wedge$ denotes
 the completion.
\end{prop}

\begin{prf}
 We may take $\mathfrak{M}_m$ as an integral model of a certain Shimura variety
 if the characteristic of $F$ is $0$ (\cite[Lemma III.4.1.1]{Harris-Taylor}), and
 as the moduli space of $\mathcal{D}$-elliptic sheaves with level structures if the characteristic of $F$ is 
 positive (\cite[Th\'eor\`eme 7.4.4]{Boyer}).
 However, there is much simpler proof; see \cite[Theorem 2.3.1]{Strauch}. Note that the generic fiber of
 the scheme $\mathfrak{M}_m$ constructed \textit{loc.\ cit.} is \'etale over
 $\Spec (\Frac V)[T_1,\ldots,T_{d-1}]$, hence purely $(d-1)$-dimensional and smooth.
\end{prf}

Next we recall the group action on the tower $(X_m)_{m\ge 0}$. Let $g$ be an element of $\GL_d(F)^{00}$ and 
$m,m'\ge 0$ integers. If we have $g^{-1}K_mg\subset K_{m'}$, then we can define a morphism
$g_{mm'}\colon X_m\longrightarrow X_{m'}$, as explained in \cite[Section 2.2]{Strauch} (under the notation
used there, it is the composite $\mathcal{M}^{(0)}_{K_m}\yrightarrow{g_{K_m}}\mathcal{M}^{(0)}_{g^{-1}K_mg}\longrightarrow \mathcal{M}^{(0)}_{K_{m'}}$). This $g_{mm'}$ is compatible with the transition maps of $(X_m)_{m\ge 0}$.
Moreover, for another $g'\in \GL_d(F)^{00}$ and $m''\ge 0$ with $g'^{-1}K_{m'}g'\subset K_{m''}$, 
we have $(gg')_{mm''}=g'_{m'm''}\circ g_{mm'}$. Therefore we get the right action of $\GL_d(F)^{00}$ on 
the tower $(X_m)_{m\ge 0}$ regarded as a pro-object (\cf \cite[Expos\'e I, 8.10]{SGA4})
of schemes over $V$.
We extend this action to $\GL_d(F)^0$ by letting the action of $\varpi I_d$ be trivial 
(here $I_d$ denotes the identity matrix).

Note that, in particular, $K_0$ acts on each $A_m$ on the left and $X_m$ on the right.
The description of this action is very simple; $g\in K_0$ maps $[(G,\iota,\phi)]\in \Def_m(A)$ to
$[(G,\iota,\phi_g)]$, where $\phi_g$ is the composite of $\phi$ and
the automorphism of $(\varpi^{-m}\mathcal{O}/\mathcal{O})^d$ induced by $g$.
Here $[-]$ denotes the isomorphism class.
Obviously $K_m$ acts on $A_m$ and $X_m$ trivially.

\begin{prop}[{{\cite[Theorem 2.1.2 (ii), Proposition 2.2.5 (i)]{Strauch}}}]
 Let $m$, $m'$ be integers with $0\le m\le m'$.
 \begin{enumerate}
  \item The morphism $X_{m',\eta}\longrightarrow X_{m,\eta}$ induced by $p_{mm'}$ is a Galois \'etale covering
	with Galois group $K_m/K_{m'}$.
  \item We have $A_{m'}^{K_m/K_{m'}}=A_m$.
 \end{enumerate}
\end{prop}

\begin{defn}\label{defn:H_LT}
 For an integer $k$, we put $H^{k,0}_{\mathrm{LT}}=\varinjlim_m (R^k\psi_{X_m}\overline{\Q}_\ell)_{x_m}$.
 The group $\GL_d(F)^0$ naturally acts on it. 
 We put $H^k_{\mathrm{LT}}=\Ind_{\GL_d(F)^0}^{\GL_d(F)}H^{k,0}_{\mathrm{LT}}$.
\end{defn}

Recall that a $\overline{\Q}_\ell$-representation $V$ of a totally disconnected locally compact group $G$ is
said to be smooth if the stabilizer of every element of $V$ is an open subgroup of $G$. 
If moreover, for every compact open subgroup $K$ of $G$ the fixed part $V^K$ is finite-dimensional,
then $V$ is said to be admissible.
An admissible $\overline{\Q}_\ell$-representation $V$ of $\GL_d(F)$ is said to be supercuspidal
if it does not appear as a subquotient of parabolically induced representations
from any proper Levi subgroup.

The main theorem of this paper is the following:

\begin{thm}\label{thm:main-thm}
 For $k\neq d-1$, the $\GL_d(F)$-representation $H^k_{\mathrm{LT}}$ 
 is admissible and has no supercuspidal representation as a subquotient.
\end{thm}

\section{Proof of the main theorem}
First we recall the natural stratification of $X_m$, which is essentially introduced in \cite{Strauch}.

\begin{defn}
 Let $m\ge 1$ and $0\le h\le d$ be integers. We denote by $\mathcal{S}_{m,h}$
 the set of $\mathcal{O}/\varpi^m\mathcal{O}$-submodules of $(\varpi^{-m}\mathcal{O}/\mathcal{O})^d$ which are direct summand of rank $h$.
 Put $\mathcal{S}_m=\bigcup_{h=0}^d \mathcal{S}_{m,h}$.

 For $J\in \mathcal{S}_{m,h}$, let $\p_J$ be the ideal of $A_m$ generated by $\phi^{\mathrm{univ}}_m(a)$
 where $a\in J$. This is known to be a prime ideal (\cite[Proposition 3.1.3]{Strauch}).
 Denote by $Y_J$ the closed subscheme of $X_m$ defined by $\p_J$. 
 Note that $Y_J\supset Y_{J'}$ for $J,J'\in \mathcal{S}_m$ with $J\subset J'$.
 Put $Z_J=Y_J\setminus \bigcup_{J'\in \mathcal{S}_m, J\subsetneq J'}Y_{J'}$, 
 $Y_m^{(h)}=\bigcup_{J\in \mathcal{S}_{m,h}}Y_J$ and $Z_m^{(h)}=\bigcup_{J\in \mathcal{S}_{m,h}}Z_J$.
 We have $Y_J=\bigcup_{J'\in \mathcal{S}_m, J\subset J'}Z_{J'}$ and
 $X_m=Y_m^{(0)}\supset Y_m^{(1)}\supset \cdots\supset Y_m^{(d)}=x_m$.
\end{defn}

\begin{lem}\label{lem:stratification}
 \begin{enumerate}
  \item We have $Y_m^{(1)}=X_{m,s}$.
  \item Let $J, J'\in \mathcal{S}_m$. If $J\neq J'$, then $Z_J\cap Z_{J'}=\varnothing$.
  \item We have $Z^{(h)}_m=Y_m^{(h)}\setminus Y_m^{(h+1)}$. In particular, $Z^{(h)}_m$ can be regarded
	as an open subscheme of $Y_m^{(h)}$.
 \end{enumerate}
\end{lem}

\begin{prf}
 \begin{enumerate}
  \item  It is sufficient to show that the ideal of $A_m$ generated by 
	 the element $\prod_{a\in (\varpi^{-m}\mathcal{O}/\mathcal{O})^d, \varpi^{m-1}a\neq 0}\phi^{\mathrm{univ}}_m(a)$ 
	 is equal to $\varpi A_m$. 
	 By the definition of a Drinfeld level structure, we have 
	 $\prod_{a\in (\varpi^{-m}\mathcal{O}/\mathcal{O})^d}(T-\phi_m^{\mathrm{univ}}(a))=[\varpi^m](T)g_m(T)$
	 for some $g_m(T)\in A_m[[T]]^\times$.
	 Since the restriction of $\phi^{\mathrm{univ}}_m$ to $(\varpi^{-m+1}\mathcal{O}/\mathcal{O})^d$
	 coincides with $\phi^{\mathrm{univ}}_{m-1}$, we have 
	 $\prod_{a\in (\varpi^{-m+1}\mathcal{O}/\mathcal{O})^d}(T-\phi_m^{\mathrm{univ}}(a))=[\varpi^{m-1}](T)g_{m-1}(T)$
	 for some $g_{m-1}(T)\in A_m[[T]]^\times$. On the other hand, since 
	 $[\varpi^m](T)=[\varpi]([\varpi^{m-1}](T))$, $\frac{[\varpi^m](T)}{[\varpi^{m-1}](T)}$ is an element of
	 $A_0[[T]]$ whose constant term is $\varpi$. Therefore, by comparing the constant terms of both sides
	 of the equation
	 \[
	 \prod_{\substack{a\in (\varpi^{-m}\mathcal{O}/\mathcal{O})^d,\\\varpi^{m-1}a\neq 0}}\bigl(T-\phi^{\mathrm{univ}}_m(a)\bigr)
	 =\frac{[\varpi^m](T)g_m(T)}{[\varpi^{m-1}](T)g_{m-1}(T)},
	 \]
	 we get the desired result.
  \item For simplicity, we write $\m$ for the maximal ideal of $A_m$.
	First note that $\phi_m^{\mathrm{univ}}\bmod \p_J\colon (\varpi^{-m}\mathcal{O}/\mathcal{O})^d\longrightarrow \m/\p_J$ is
	$0$ on $J$. 

	We may assume that $\rank J\le \rank J'$.
	Take a primary decomposition $\p_J+\p_{J'}=\bigcap_{i=1}^l\mathfrak{q}_i$, where $\mathfrak{q}_i$ is
	a primary ideal of $A_m$. Put $\p_i=\sqrt{\mathfrak{q}_i}$. Then $\p_i$ is a prime ideal and
	we have $\sqrt{\p_J+\p_{J'}}=\bigcap_{i=1}^l\p_i$.
	Let $J''_i$ be the kernel of 
	$\phi_m^{\mathrm{univ}}\bmod \p_i\colon (\varpi^{-m}\mathcal{O}/\mathcal{O})^d\longrightarrow \m/\p_i$.
	By the subsequent lemma, it is an element of $\mathcal{S}_m$.
	Moreover, since $\p_J, \p_{J'}\subset \p_i$, $J''_i$ contains $J$ and $J'$.
	Therefore $J\subsetneq J''_i$, for $\rank J\le \rank J'$ and $J\neq J'$.
	By the definition we have $\p_{J''_i}\subset \p_i$. Thus we have 
	$\p_{J''_1}\cdots \p_{J''_l}\subset \bigcap_{i=1}^l \p_i=\sqrt{\p_J+\p_{J'}}$,
	which implies the set-theoretical inclusion
	$Y_J\cap Y_{J'}\subset \bigcup_{i=1}^l Y_{J''_i}\subset Y_J\setminus Z_J$.
	Hence we have $Z_J\cap Z_{J'}=\varnothing$, as desired.
  \item It is an easy consequence of ii) and $Y_J=\bigcup_{J'\in \mathcal{S}_m, J\subset J'}Z_{J'}$. 
 \end{enumerate}
\end{prf}

\begin{lem}\label{lem:kernel-of-Drinfeld-level-str}
 Let $A$ be an object of the category $\mathcal{C}$, $m\ge 1$ an integer, 
 and $(G,\iota,\phi)$ a triple that appeared in the definition of $\Def_m(A)$. 
 Assume that $A$ is an integral domain. Then $\Ker \phi$ is a direct summand of $(\varpi^{-m}\mathcal{O}/\mathcal{O})^d$.
\end{lem}

\begin{prf}
 First note that a finitely generated $\mathcal{O}/\varpi^{m}\mathcal{O}$-module $M$ is free if and only if the following
 condition holds:
 \begin{quote}
  for every $x\in M$ with $\varpi x=0$, there exists an element $y\in M$ such that $x=\varpi^{m-1} y$.
 \end{quote}
 Indeed, it is easy to see that the condition above is equivalent to $\Tor_1(M,\mathcal{O}/\varpi\mathcal{O})=0$, which is
 equivalent to the flatness of $M$.

 In order to show that $\Ker \phi$ is a direct summand, it suffices to show that $\Imm \phi$ is a free
 $\mathcal{O}/\varpi^m\mathcal{O}$-module. We will verify the condition above. Assume that $[\varpi](\phi(a))=0$.
 Since $\prod_{b\in (\varpi^{-1}\mathcal{O}/\mathcal{O})^d}(T-\phi(b))=[\varpi](T)g(T)$ for some $g(T)\in A[[T]]^\times$,
 we have $\prod_{b\in (\varpi^{-1}\mathcal{O}/\mathcal{O})^d}(\phi(a)-\phi(b))=0$.
 Since $A$ is an integral domain, there exists an element $b\in (\varpi^{-1}\mathcal{O}/\mathcal{O})^d$ such that
 $\phi(a)=\phi(b)$. Take $b'\in (\varpi^{-m}\mathcal{O}/\mathcal{O})^d$ such that $\varpi^{m-1}b'=b$. 
 Then $[\varpi^{m-1}](\phi(b'))=\phi(\varpi^{m-1}b')=\phi(b)=\phi(a)$. This completes the proof.
\end{prf}

\begin{prop}\label{prop:finiteness-LT}
 For every pair of integers $k$ and $h\ge 1$, $H^k_{Z_m^{(h)}}(X_{m,s},R\psi_{X_m}\overline{\Q}_\ell)$
 is a finite-dimensional $\overline{\Q}_\ell$-vector space.
\end{prop}

\begin{prf}
 Clear from Proposition \ref{prop:finiteness-local-coh} and Proposition \ref{prop:algebraization-A_m}.
\end{prf}

\begin{rem}
 It is possible to take an algebraization $\mathfrak{M}_m$ in Proposition \ref{prop:algebraization-A_m}
 so that the locally closed subscheme
 $Z_m^{(h)}$ comes from a locally closed subscheme of $\mathfrak{M}_m$. Then we may use 
 Remark \ref{rem:without-Gabber} in place of Proposition \ref{prop:finiteness-local-coh}.
\end{rem}

As explained in \cite[Section 3.1]{Strauch}, this stratification is preserved by the action of $\GL_d(F)^0$
in the sense of the following.

Fix an integer $h$ with $1\le h\le d-1$. Denote by $J_h^0$ (resp.\ $J_{m,h}^0$ for $m\ge 1$)
the direct summand of $\mathcal{O}^d$ (resp.\ $(\varpi^{-m}\mathcal{O}/\mathcal{O})^d$) generated by the first $h$ elements of
the standard basis. Let $P_h(\mathcal{O})\subset \GL_d(\mathcal{O})$ be the stabilizer of $J_h^0$. Similarly we define
$P_h(F)\subset \GL_d(F)$. Put $P_h(F)^0=P_h(F)\cap \GL_d(F)^0$.

Obviously $\GL_d(\mathcal{O})$ acts on $\mathcal{S}_{m,h}$. By mapping $g$ to $gJ_{m,h}^0$, we have a map
$\GL_d(\mathcal{O})\longrightarrow \mathcal{S}_{m,h}$, which induces natural isomorphisms
\[
 \mathcal{S}_{m,h}\cong K_m\backslash \GL_d(\mathcal{O})/P_h(\mathcal{O})\cong K_m\backslash \GL_d(F)/P_h(F)
 \cong K_m\backslash \GL_d(F)^0/P_h(F)^0.
\]
Moreover, for $1\le m\le m'$, it is easy to see that under the isomorphism above the projection 
$K_{m'}\backslash \GL_d(F)^0/P_h(F)^0\longrightarrow K_m\backslash \GL_d(F)^0/P_h(F)^0$
corresponds to the map $\mathcal{S}_{m',h}\longrightarrow \mathcal{S}_{m,h}$ that maps $J$ to
$J\cap (\varpi^{-m}\mathcal{O}/\mathcal{O})^d$.

Let $g$ be an element of $\GL_d(F)^0$ and $m$, $m'$ positive integers with $g^{-1}K_mg\subset K_{m'}$.
Then we can define the map $\mathcal{S}_{m,h}\longrightarrow \mathcal{S}_{m',h}$ such that
the corresponding map $K_m\backslash \GL_d(F)^0/P_h(F)^0\longrightarrow K_{m'}\backslash \GL_d(F)^0/P_h(F)^0$
is given by $K_mxP_h(F)^0\longmapsto K_{m'}g^{-1}xP_h(F)^0$. 
We denote the image of $J\in \mathcal{S}_{m,h}$ under this map by $g^{-1}J$. Then, the morphism
$g_{mm'}\colon X_m\longrightarrow X_{m'}$ maps $Z_J$ into $Z_{g^{-1}J}$.

\begin{defn}
 Let $h$ be an integer with $1\le h\le d-1$. For integers $m$, $m'$ with $1\le m\le m'$, 
 $Z_{J^0_{m',h}}$ is open and closed in $p_{mm'}^{-1}(Z_{J^0_{m,h}})$, and thus we have the inductive system
 $(H^k_{Z_{J^0_{m,h}}}(X_{m,s},R\psi_{X_m}\overline{\Q}_\ell))_{m\ge 1}$.
 Put $W^{(h,k,0)}=\varinjlim_m H^k_{Z_{J^0_{m,h}}}(X_{m,s},R\psi_{X_m}\overline{\Q}_\ell)$.
 Since the image of $H^k_{Z_{J^0_{m,h}}}(X_{m,s},R\psi_{X_m}\overline{\Q}_\ell)$ in $W^{(h,k,0)}$ is
 fixed by the open subgroup $P_h(F)^0\cap K_m$ of $P_h(F)^0$, 
 $W^{(h,k,0)}$ is a smooth representation of $P_h(F)^0$.
 Similarly we may define $V^{(h,k,0)}$ as $\varinjlim_m H^k_{Z_m^{(h)}}(X_{m,s},R\psi_{X_m}\overline{\Q}_\ell)$,
 which is a smooth representation of $\GL_d(F)^0$.
 Put $W^{(h,k)}=\Ind_{P_h(F)^0}^{P_h(F)} W^{(h,k,0)}$ and $V^{(h,k)}=\Ind_{\GL_d(F)^0}^{\GL_d(F)}V^{(h,k,0)}$.
\end{defn}

\begin{prop}\label{prop:V-W}
 Let $k$ and $h$ be integers with $1\le h\le d-1$.
 \begin{enumerate}
  \item The $P_h(F)^0$-representation $W^{(h,k,0)}$ is admissible.
  \item We have an isomorphism $V^{(h,k,0)}\cong \Ind_{P_h(F)^0}^{\GL_d(F)^0}W^{(h,k,0)}$.
  \item The unipotent radical of $P_h(F)$ acts on $W^{(h,k)}$ trivially.
  \item There is no supercuspidal subquotient of $V^{(h,k)}$.
 \end{enumerate}
\end{prop}

\begin{prf}
 \begin{enumerate}
  \item For an integer $m\ge 1$ put $P_h(F)_m^0=P_h(F)^0\cap K_m$.
	It is sufficient to show that $(W^{(h,k,0)})^{P_h(F)^0_m}$ is 
	a finite-dimensional $\overline{\Q}_\ell$-vector space for every $m$.
	Let $m'$ be an integer with $m'\ge m$.
	Since $p_{mm'}^{-1}(Z_{J_{m,h}^0})$ is the disjoint union of $(Z_{J_{m',h}^0})g^{-1}=Z_{gJ_{m',h}^0}$
	where $g$ runs through $K_m/K_{m'}$ (Lemma \ref{lem:stratification} ii)), we have the natural isomorphisms
	\begin{align*}
	 H^k_{p_{mm'}^{-1}(Z_{J^0_{m,h}})}(X_{m',s},R\psi_{X_{m'}}\overline{\Q}_\ell)
	 &\cong \bigoplus_{g\in K_m/K_{m'}}
	 H^k_{gZ_{J^0_{m',h}}}(X_{m',s},R\psi_{X_{m'}}\overline{\Q}_\ell)\\
	 &\cong \Ind_{P_h(F)^0_m/P_h(F)^0_{m'}}^{K_m/K_{m'}}H^k_{Z_{J^0_{m',h}}}(X_{m',s},R\psi_{X_{m'}}\overline{\Q}_\ell)
	\end{align*}
	(the second isomorphism is the one given in the proof of \cite[Lemme 13.2]{Boyer}).
	Since  $K_m/K_{m'}$ is a $p$-group where $p$ is the characteristic of $\F_q$
	and $A_m=A_{m'}^{K_m/K_{m'}}$, by Proposition \ref{prop:invariant} we have
	\begin{align*}
	 &H^k_{Z_{J^0_{m',h}}}(X_{m',s},R\psi_{X_{m'}}\overline{\Q}_\ell)^{P_h(F)^0_m}\cong
	\Bigl(\Ind_{P_h(F)^0_m/P_h(F)^0_{m'}}^{K_m/K_{m'}}H^k_{Z_{J^0_{m',h}}}(X_{m',s},R\psi_{X_{m'}}\overline{\Q}_\ell)\Bigr)^{K_m/K_{m'}}\\
	 &\qquad \cong  H^k_{p_{mm'}^{-1}(Z_{J^0_{m,h}})}(X_{m',s},R\psi_{X_{m'}}\overline{\Q}_\ell)^{K_m/K_{m'}}
	 \cong H^k_{Z_{J^0_{m,h}}}(X_{m,s},R\psi_{X_m}\overline{\Q}_\ell).
	\end{align*}
	Moreover it is not difficult to see that this isomorphism is induced by the natural pull-back map
	$H^k_{Z_{J^0_{m,h}}}(X_{m,s},R\psi_{X_m}\overline{\Q}_\ell)\yrightarrow{p_{mm'}^*} H^k_{Z_{J^0_{m',h}}}(X_{m',s},R\psi_{X_{m'}}\overline{\Q}_\ell)$. Thus by taking the inductive limit with respect to $m'$, we have 
	\[
	(W^{(h,k,0)})^{P_h(F)^0_m}\cong H^k_{Z_{J^0_{m,h}}}(X_{m,s},R\psi_{X_m}\overline{\Q}_\ell).
	\]
	Since the right hand side is a finite-dimensional $\overline{\Q}_\ell$-vector space
	by Proposition \ref{prop:finiteness-LT}, so is the left hand side.
  \item It can be proved in the same way as \cite[Theorem 4.3.2 (ii)]{Strauch}.
  \item By i), $W^{(h,k)}$ is an admissible $P_h(F)$-representation. Thus by \cite[Lemme 13.2.3]{Boyer}, 
	its restriction to the unipotent radical of $P_h(F)$ is trivial.
  \item By ii), we have $V^{(h,k)}\cong \Ind_{P_h(F)}^{\GL_d(F)}W^{(h,k)}$. By iii), the right hand side
	is obtained by a parabolic induction from the Levi subgroup of $P_h(F)$. 
	Thus it has no supercuspidal subquotient.
 \end{enumerate}
\end{prf}

Now we can give a proof of our main theorem.

\begin{prfwithname}[of Theorem \ref{thm:main-thm}]
 For integers $k$ and $h$ with $1\le h\le d$,
 we put $V'^{(h,k,0)}=\varinjlim_m H^k_{Y_m^{(h)}}(X_{m,s},R\psi_{X_m}\overline{\Q}_\ell)$ and
 $V'^{(h,k)}=\Ind_{\GL_d(F)^0}^{\GL_d(F)}V'^{(h,k,0)}$ (since $p_{mm'}^{-1}(Y_m^{(h)})=Y_{m'}^{(h)}$ for
 $1\le m\le m'$, we have the inductive system $(H^k_{Y_m^{(h)}}(X_{m,s},R\psi_{X_m}\overline{\Q}_\ell))_{m\ge 1}$). It is easy to see that 
 $V'^{(h,k,0)}$ (resp.\ $V'^{(h,k)}$) is a smooth representation of $\GL_d(F)^0$ (resp.\ $\GL_d(F)$).
 For an integer $h$ with $0\le h\le d-1$, the exact sequence
 \begin{align*}
  \cdots&\longrightarrow H^k_{Y_m^{(h+1)}}(X_{m,s},R\psi_{X_m}\overline{\Q}_\ell)
  \longrightarrow H^k_{Y_m^{(h)}}(X_{m,s},R\psi_{X_m}\overline{\Q}_\ell)\\
  &\qquad\longrightarrow H^k_{Z_m^{(h)}}(X_{m,s},R\psi_{X_m}\overline{\Q}_\ell)
  \longrightarrow H^{k+1}_{Y_m^{(h+1)}}(X_{m,s},R\psi_{X_m}\overline{\Q}_\ell)\longrightarrow \cdots
 \end{align*}
 gives the exact sequence $V'^{(h+1,k)}\longrightarrow V'^{(h,k)}\longrightarrow V^{(h,k)}$.

 Assume that $k<d-1$. Then by Proposition \ref{prop:coh-dim} and Proposition \ref{prop:algebraization-A_m}, 
 \[
  V'^{(d,k)}=\Ind_{\GL_d(F)^0}^{\GL_d(F)}\varinjlim_m H^k_{x_m}(X_{m,s},R\psi_{X_m}\overline{\Q}_\ell)=0.
 \]
 Therefore, by Proposition \ref{prop:V-W} iv) and the exact sequence above, we have inductively that
 for every $h$ with $1\le h\le d$, $V'^{(h,k)}$ is an admissible $\GL_d(F)$-representation which has
 no supercuspidal subquotient. On the other hand, by Lemma \ref{lem:stratification} i),
 \[
  V'^{(1,k)}=\Ind_{\GL_d(F)^0}^{\GL_d(F)}\varinjlim_m H^k(X_{m,s},R\psi_{X_m}\overline{\Q}_\ell)=
 \Ind_{\GL_d(F)^0}^{\GL_d(F)}\varinjlim_m(R^k\psi_{X_m}\overline{\Q}_\ell)_{x_m}=H^k_{\mathrm{LT}}
 \]
 (note that $X_{m,s}$ is the spectrum of a strict local ring with closed point $x_m$). 
 Hence $H^k_{\mathrm{LT}}$ is an admissible $\GL_d(F)$-representation which has
 no supercuspidal subquotient.

 Next assume that $k>d-1$. Then by Proposition \ref{prop:coh-dim} and Proposition \ref{prop:algebraization-A_m}, 
 $(R^k\psi_{X_m}\overline{\Q}_\ell)_{x_m}=0$.
 Therefore $H^k_{\mathrm{LT}}=0$ and Theorem \ref{thm:main-thm} holds obviously.
\end{prfwithname}

\end{document}